\documentclass{amsart}

\usepackage{amsmath,amssymb}

\newtheorem{lem}{Lemma}
\newtheorem{cor}{Corollary}
\newtheorem{thm}{Theorem}
\newtheorem{rem}{Remark}

\def\Cov#1#2{{\rm Cov}_{#1}(#2)}
\def\Cap#1#2{{\rm Cap}_{#1}(#2)}

\def\vol{{\rm vol}}
\def\dist{{\rm dist}}

\title{SEPARATION THEOREMS FOR COMPACT HAUSDORFF FOLIATIONS}
\author{Wojciech Koz\l owski \and Szymon M. Walczak}
\date{\today}

\begin{document}
\baselineskip=16pt

\begin{abstract}
We investigate compact Hausdorff foliations on compact Riemannian manifolds in the context of the Gromov-Hausdorff distance theory. We give some sufficient conditions for such foliations to be separated in the Gromov-Hausdorff topology (GH-separation theorem).
\end{abstract}

\maketitle

\section{Introduction}

The concept of the Gromov-Hausdorff distance (briefly GH-distance), being a generalization of the notion of the Hausdorff distance, was originally introduced by M. Gromov \cite{G} in the late  1970s. Next,  Gromov, Katsuda, Peters and others showed that the GH-distance theory applied to the Riemannian manifolds leads to remarkable results \cite{CY}, e.g., the {\em Cheeger's finiteness theorem} follows from the {\em Gromov-Katsuda convergence theorem} \cite{Ka}. 

In \cite{W} and \cite{W1}, the second named author investigated warped compact Hausdorff foliations from the GH-distance theory point of view. He gave a necessary and sufficient conditions for a sequence of warped compact Hausdorff foliation to be converged to the space of leaves with quotient metric. 

In the light of the results there appears a natural question: \textit{Suppose that two compact metric spaces $X$ and $X'$ are GH-close. Are always the compact Hausdorff foliations $(M,\mathcal{F},g)$ and $(M',\mathcal{F}',g)'$ with space of leaves coinciding with $X$ and $X'$, respectively, GH-close?}

In this paper we show that for compact Hausdorff foliations the answer is negative ({\em GH-separation Theorem} - Theorem \ref{CompactHausdorffSeparation} - the main result of the paper).


\section{Gromov-Hausdorff distance}

Let $C,K\subset X$ be compact subsets of a metric space $(X,d)$. The number 
\begin{equation*}
d_{{\rm H}} (C,K)=\inf \{\epsilon >0: C\subset N(K,\epsilon) \wedge K\subset N(C,\epsilon)\}, 
\end{equation*}
where $N(A,\epsilon) = \{x\in X: d(x,A)<\epsilon\}$, is called the {\em Hausdorff distance} of $C$ and $K$. 

Let $(X, d_X)$ and $(Y, d_Y)$ be arbitrary compact metric spaces. Equip the disjoint union $X\amalg Y$ with an admissible metric $d$, i.e., the metric which extends $d_X$ and $d_Y$. The {\em Gromov-Hausdorff distance} (cf.  \cite{CY}, \cite{G} and \cite{P}) (GH-distance) $d_{{\rm GH}}$ of the spaces $(X, d_X)$ and $(Y, d_Y)$ one can define as
\begin{equation*}
d_{{\rm GH}} (X,Y):={\rm inf}\{d_{\rm H}(X,Y)\},
\end{equation*}
where the infimum is taken over all admissible metrics on $X\amalg Y$. Note that two compact metric spaces are isometric iff their GH-distance equals zero. Consequently, GH-distance in the class of all classes of isometry of compact metric spaces with the GH-distance is a metric.

\begin{lem}\label{GHlemma}
If there exist $\epsilon$-nets $\{x_1,\dots ,x_k\}\subset X$ and $\{y_1,\dots ,y_k\}\subset Y$ satisfying for all $1\leq i,j\leq k$ 
\begin{equation*}
|d_X(x_i,x_j) - d_Y(y_i,y_j)|\leq \epsilon
\end{equation*}
then $d_{{\rm GH}} (X,Y)\leq 3\epsilon$.
\end{lem}

\begin{proof}
For a proof we refer to \cite{CY}.
\end{proof}

\begin{lem}\label{GHlemma-converse}
If $d_{{\rm GH}} (X,Y) \leq \epsilon$ then for every $\epsilon$-net $\{x_1,\dots ,x_k\}\subset X$ there 
exists a $3\epsilon$-net $\{y_1,\dots ,y_k\}\subset Y$ such that $|d_X(x_i,x_j) - d_Y(y_i,y_j)|\leq 2\epsilon$ for all $1\leq i,j\leq k$.
\end{lem}
\begin{proof}
Since $d_{{\rm GH}}(X,Y) \leq \epsilon$ there exists an admissible metric $d$ on $X\amalg Y$ such that 
the Hausdorff distance $d_{{\rm GH}}(X,Y)\leq d_{\rm H} (X,Y) \leq \epsilon$. Let $\{x_1,\dots ,x_k\}$ 
be an $\epsilon$-net on $X$. For every $i\in\{1,\dots,k\}$ there exists $y_i\in Y$, $y_i\neq y_j$ while 
$i\neq j$, such that $d(x_i,y_i) \leq \epsilon$. Since $d$ is an extension of metrics $d_X$ and $d_Y$, $\{y_1,\dots,y_k\}$ is a $3\epsilon$-net on $Y$. Moreover,
\begin{eqnarray*}
d_X(x_i,x_j) \leq d(x_i,y_i) + d(y_i,y_j) + d(y_j,x_j) \leq  d_Y(y_i,y_j) + 2\epsilon,
\end{eqnarray*}
and similarly $d_Y(y_i,y_j)\leq  d_X(x_i,x_j)+2\epsilon$.
\end{proof}

Let $\Cov{\epsilon}{X}$ denote the smallest number of open $\epsilon$-balls which covers $X$, and $\Cap{\epsilon}{X}$ the largest number of disjoint $\epsilon$-balls contained in $X$. Obviously  
\begin{equation*}
\Cap{\epsilon}{X}\leq\Cov{\epsilon}{X}. 
\end{equation*}

\begin{lem}\label{CorCoveringTrick}
Let $(X,d)$ be a compact metric space. Then $\Cov{2\delta}{M} \leq \Cap{\delta}{M}$. More precisely, if $x_1,\dots,x_ {\Cap{\delta}{M}}$ are the centres of disjoint balls $B(x_i,r)$ then the balls $B(x_i,2r)$ cover $M$. 
\end{lem}
\begin{proof}
Let $k=\Cap{\delta}{X}$ and $B(x_1,\delta),\dots,B(x_1,\delta)$ be a family of open disjoint balls in $X$. Let $x\in X$. Then $B(x,\delta)\cap B(x_i,\delta)\neq \emptyset$, and $d(x,y) < \delta$ and $d(y,x_i)<\delta$. Therefore, $x\in  B(x_i, 2\delta)$ for some $i\in\{1,\dots,k\}$. Hence, $B(x_1,2\delta),\dots,B(x_1,2\delta)$ cover $X$.
\end{proof}

\begin{lem}\label{LemmaCapCov}
If $\Cov{\epsilon}{X} < \Cap{3\epsilon}{Y}$ then $d_{{\rm GH}}(X,Y)>\epsilon$.
\end{lem}
\begin{proof}
Suppose that $d_{{\rm GH}}(X,Y)\leq \epsilon$. Let $\{x_1,\dots,x_k\}$, $k = \Cov{\epsilon}{X}$, be an $\epsilon$-net on $X$. Then, by Lemma \ref{GHlemma-converse}, there exists a $3\epsilon$-net $\{y_1,\dots ,y_k\}$ on $Y$, so $\Cov{3\epsilon}{Y}\leq \Cov{\epsilon}{X}$. Thus 
\begin{equation*}
\Cov{3\epsilon}{Y}\leq \Cov{\epsilon}{X} < \Cap{3\epsilon}{Y} \leq \Cov{3\epsilon}{Y}. 
\end{equation*}
Contradiction gives us the statement.
\end{proof}

Let $p\geq 0$. A Borel measure $\mu$ on a metric space $(X,d)$ is called \textit{a $p$-dimensional Bishop measure on $(X,d)$} if there exist constants $\beta\geq 1$ and $\eta_0>0$ such that for all $\eta<\eta_0$ and every $x\in X$
\begin{equation}\label{BishopInequality}
\frac{1}{\beta} \eta^p\leq \mu(B(x,\eta)) \leq \beta\eta^p,
\end{equation}
where $B(x,\eta) = \{y\in X: d(x,y)<\eta\}$. 

Let $(X,d)$ be a length-space, i.e. $d(x,y) = \inf \{l(\gamma)\}$, where $\gamma:[0,1]\to X$ is a curve such that $\gamma(0)=x$, $\gamma(1)=y$, and $l(\gamma)$ denotes the length of $\gamma$.

\begin{lem}\label{Lemma2} 
If the balls $B(x,\delta)$ and $B(y,\delta)$ are disjoint, then $d(x,y)\geq 2\delta$.
\end{lem}
\begin{proof} Suppose that $d(x,y)< 2\delta$. Let $\gamma:[0,1]\to X$ be a curve from $x=\gamma(0)$ to $y=\gamma(1)$ with its length $l(\gamma)<2\delta$. Let $t_0\in[0,1]$ be such that
\begin{equation*}
l(\gamma|{[0,t_0]}) = l(\gamma|{[t_0, 1]}) = \frac{1}{2}l(\gamma) \leq \delta. 
\end{equation*}
If $z=\gamma(1/2)$ then $d(x,z) < \delta$ and $d(z,y) < \delta$. Thus, $z\in B(x,\delta)\cap B(y, \delta)$. Contradiction ends our proof.
\end{proof}

\begin{lem}\label{LemmaCapEstimation}
Let $(X,d)$ be a compact length space, $p\geq 1$, and let $\mu$ be a $p$-dimensional Bishop measure on $(X,d)$ with constants $\beta>1$, $\eta_0>0$. There exist positive constants $C\geq 1$ and $\theta>0$ such that for every $0<r<\theta$ and $x\in M$,
\begin{equation*}
 \frac{1}{C r^p} \mu(X) \leq \Cap{r}{X} \leq \frac{C}{r^p} \mu(X).
\end{equation*}
\end{lem}
\begin{proof} Let $0<r<\eta_0$, $k=\Cap{r}{X}$, and let $B(x_1,\delta),\dots,B(x_k,\delta)$ be a family of open disjoint balls in $X$. By (\ref{BishopInequality}), $\mu (B(x_i,r))\geq \beta^{-1} r^p $, and
\begin{equation}\label{KeyLemmaEq1}
\mu(X) \geq \sum_{i=1}^k \mu(B(x_i,r)) \geq k\cdot \frac{1}{\beta} r^p.
\end{equation}
Let $0<r<\eta_0/2$. By Lemma \ref{CorCoveringTrick}, we have
\begin{equation}\label{KeyLemmaEq2}
\mu(X) \leq \sum_{i=1}^k \mu(B(x_i,2r)) \leq k\cdot \beta (2r)^p .
\end{equation}
Putting $C=\beta2^p$ and $\theta=\eta_0/2$, (\ref{KeyLemmaEq1}) and (\ref{KeyLemmaEq2}) give us the statement.
\end{proof}

\begin{cor}\label{CorollaryCapEstimation}
Let $0<r<\theta$ and $\alpha>0$ be such that $r\alpha<\theta$. Then
\begin{equation*}
\alpha^{-p}C^{-2} \Cap{r}{X} \leq \Cap{\alpha r}{X} \leq \alpha^{-p}C^{2} \Cap{r}{X}. 
\end{equation*}
\end{cor}

\begin{rem}
Note that the volume form on a $n$-dimensional compact Riemannian manifold defines an $n$-dimensional Bishop measure.
\end{rem}

\section{Compact Hausdorff foliations}

A foliation with all leaves compact is called a \textit{compact foliation}. Let us consider any compact foliation $\mathcal{F}$ on a manifold $M$, and let $\pi:M\to \mathcal{L}$ denote a quotient map onto the space of leaves $\mathcal{L}$, this means that $\pi$ identifies each leaf to a point. The space of leaves often is non-Hausdorff. Due to the results by D.B.A. Epstein \cite{E2}, we recall theorems that describe the topology of such foliation: 
\begin{thm}\label{epsteinconditions1}
The following conditions are equivalent.
\begin{itemize}
\item[(i)] $\pi$ is a closed map.
\item[(ii)] $\pi$ maps compact sets onto closed sets.
\item[(iii)] Each leaf has arbitrarily small saturated neighbourhoods.
\item[(iv)] $\mathcal{L}$ with quotient topology is Hausdorff.
\item[(v)] If $K\subseteq M$ is compact, then the saturation of $K$ is also compact.
\end{itemize}
\end{thm}
\begin{proof}
For a proof we refer to \cite{E2}, Theorem 4.1.
\end{proof}

Let $M$ be a Riemannian manifold and $N$ a submanifold on $M$. One can consider the induced Riemannian structure on $N$ and introduce a volume of $N$ as it's volume $\vol_{N}$ in the induced Riemannian structure. The next theorem describes the relation between the volume of the leaves defined above (briefly the volume function), the holonomy group of a leaf, and the topology of the space of leaves of a foliation $\mathcal{F}$ on a Riemannian manifold $(M,g)$.
\begin{thm}\label{epsteinconditions2}
If $(M,\mathcal{F},g)$ is a foliated Riemannian manifold and $L$ is a compact leaf of $\mathcal{F}$, then the following conditions are equivalent.
\begin{itemize}
\item[(i)] There exists a saturated neighbourhood $N$ of the leaf $L$ such that the volume function is bounded on $N$.
\item[(ii)] The holonomy group of $L$ is finite.
\end{itemize}
\end{thm}
\begin{proof}
For a proof we refer to \cite{E2}, [Theorem 4.2].
\end{proof}

The conditions of Theorem \ref{epsteinconditions2} imply that some saturated neighbourhood $U$ of a compact leaf $L$ consists of compact leaves, and in $U$ the conditions of Theorem \ref{epsteinconditions1} are satisfied in $U$. Moreover, by Reeb Stability Theorem, on a foliated manifold the conditions of Theorem \ref{epsteinconditions1} imply the conditions of Theorem \ref{epsteinconditions2}.

A compact foliation which space of leaves is Hausdorff is called \textit{compact Hausdorff foliation}. As an easy corollary of the above theorems we have:
\begin{cor}
Let $(M,\mathcal{F},g)$ be a compact Riemannian manifold carrying compact Hausdorff foliation. Then $\sup\limits_{L\in\mathcal{F}} \vol(L) < \infty $.
\end{cor}

Now, let us consider the space of leaves $\mathcal{L}$ of an arbitrary compact Hausdorff foliation on a compact Riemannian manifold. Let us introduce on $\mathcal{L}$ a metric $\rho$ defined by
\begin{equation*}\label{rhodef}
\rho(L,L') = \inf \{\sum_{i=1}^{n-1} \dist (L_i,L_{i+1})\},
\end{equation*}
where $L_1=L$, $L_n=L'$, and the infimum is taken over all finite sequences of leaves. One can see that for a compact Riemannian foliation $\mathcal{F}$ the distance $\rho$ coincides with Hausdorff distance of leaves of $\mathcal{F}$. 

\begin{rem}\label{rhoRemark}
Let $g,g'$ be two Riemannian metrics on a compact foliated manifold $(M,mathcal{F})$, where $\mathcal{F}$ is a compact Hausdorff foliation. Denote by $\rho$ and $\rho'$ two metrics on the space of leaves constructed using $g$ and $g'$, respectively. Since $M$ is compact, then $\frac{1}{C} g \leq g' \leq Cg$ for some constant $C\geq 1$. One can check that 
\[
\frac{1}{C} \rho \leq \rho'\leq C \rho.
\]
\end{rem}

In further considerations we will need the following:

\begin{lem}\label{RiemannianFoliation}
For every compact Hausdorff foliation $\mathcal{F}$ on a compact Riemannian manifold $(M,g)$ there exists Riemannian structure $\tilde g$ on $M$ such that $\mathcal{F}$ becomes a Riemannian foliation, and for any leaf $L\in\mathcal{F}$ we have $\tilde g|L = g|L$ 
\end{lem}
\begin{proof}
Obvious. See \cite{MM}.
\end{proof}

\section{Separation Theorem}

We say that a compact metric space $(X',d')$ is \textit{broader} than a metric space $(X,d)$, and we briefly write $X'\succeq X$,  if $\Cap{\delta}{X'}\geq \Cap{\delta}{X}$ for all $\delta>0$.

Let $d>0$ be a real number. Let us denote by $\mathcal{M}(d,C,p,n)$ the class of all $n$-dimensional compact foliated Riemannian manifolds $(M,\mathcal{F},g)$ carrying a compact Hausdorff foliation of dimension $p$ satisfying:
\begin{enumerate}
\item For any leaf $L\in\mathcal{F}$, $\epsilon<d$, and any two balls $B_L(x,\epsilon),B_L(y,\epsilon)$ that are disjoint in $L$, the balls $B(x,\epsilon)$ and $B(y,\epsilon)$ are disjoint in $M$;
\item $\max\limits_{L\in\mathcal{F}} C_L \leq C$, $\min\limits_{L\in\mathcal{F}} \theta_L \leq C$, $C_M\leq C$, $\theta_M\leq C$, where $C_L,\theta_L, C_M, \theta_M$ are the constants mentioned in Lemma \ref{LemmaCapEstimation} for a leaf $L$ of $\mathcal{F}$ and for the manifold $M$, respectively;
\item $\frac{1}{C}\leq \vol(L)\leq C$ for all $L\in\mathcal{F}$;
\item There exists a Riemannian structure $\tilde g$ on $M$ satisfying $\frac{1}{C} g \leq \tilde g \leq Cg$ such that on $(M,\tilde g)$ the foliation $\mathcal{F}$ becomes a compact Riemannian foliation.
\end{enumerate}

Let $d>0$, $C\geq 1$, and let $p,p',n,n'\in \mathbb{N}$ be such that $p'>p$ and $n'\geq n$.

\begin{thm}\label{CompactHausdorffSeparation} {\rm [GH-separation Theorem]} 
There exists $\epsilon>0$ such that for any $(M,\mathcal{F},g)\in \mathcal{M}(d_0,C,p,n)$ and $(M',\mathcal{F}',g')\in \mathcal{M}(d_0,C,p',n')$ such that $(M'/\mathcal{F}',\rho')\succeq (M/\mathcal{F},\rho)$ we have $d_{{\rm GH}}(M,M')>\epsilon$.
\end{thm}
\begin{proof}
Let $\mathcal{L}=M/\mathcal{F}$, $\mathcal{L}'=M'/\mathcal{F'}$, and let $\pi:M\to \mathcal{L}$ and $\pi':M'\to \mathcal{L}'$ denote the natural projections. Let $r<\min\{d,C\}/3C$, and let $B(x_1,r/2),\dots,B(x_k,r/2)$, $k=\Cap{r/2}{\mathcal{L}}$, be a family of open disjoint balls in $\mathcal{L}$. Since $\mathcal{L}'\succeq \mathcal{L}$ we can choose points $x_1',\dots,x_k'$ in $\mathcal{L}'$ such that the balls $B(x'_1,r/2),\dots,B(x'_k,r/2)$ are also disjoint.

Now, in every leaf $L_i'=(\pi')^{-1}(x'_i)$ let us choose points 
\[
\xi_{i,1}',\dots,\xi_{i,l'}',
\]
where $l'=\min\limits_{L\in\mathcal{F}'} \Cap{r/2}{L}>0$, such that the balls $B_{L'_i}(x'_1,r/2),\dots,B_{L'_i}(x'_l,r/2)$ are disjoint. Since $r<d$, the balls $B(\xi_{i,j}',r/2)$ are disjoint in $M'$. Consequently, by Lemma \ref{LemmaCapEstimation}, we have
\begin{equation}\label{AConstant}
\Cap{\frac{r}{2}}{M'}\geq k\cdot l'\geq k\cdot \frac{2^{p'}}{\max\limits_{L\in\mathcal{F}'} C_{L} r^{p'}} \min\limits_{L\in\mathcal{F}'}\vol'(L) \geq  \frac{k}{C^2}\cdot\frac{2^{p'}}{r^{p'}},
\end{equation}
where $\vol'(L)$ denote the volume of a leaf in the induced Reimannian structure.

Now, let $\tilde g$ be a Riemannian structure on $M$ mentioned in Lemma \ref{RiemannianFoliation} such that $\mathcal{F}$ becomes a Riemannian foliation and such that
\begin{equation}\label{tildeConstant}
\frac{1}{C} g \leq \tilde g \leq Cg.
\end{equation}
Let us choose in $L_i=\pi^{-1}(x_i)$ points $\xi_{i,1},\dots,\xi_{i,l_i}$, $l_i=\Cap{r/2}{L_i}$ such the balls $B_{L_i}(\xi_{i,j},r/2)$ are pairewise disjoint on $(L_i,\tilde g|_{L_i})$. By (\ref{tildeConstant}) and Remark \ref{rhoRemark}, the balls $B(\xi_{i,j},Cr)$ covers $(M,g)$, and
\begin{equation*}
\Cov{Cr}{M}\leq k\cdot l,
\end{equation*}
where $l=\max\limits_{i\in\{1,\dots,k\}}\{l_i\}$. Moreover, by (\ref{tildeConstant}) and Lemma \ref{LemmaCapEstimation},
\begin{equation*}
l \leq \max\limits_{L\in\mathcal{F}}\Cap{\frac{r}{2C}}{L} \leq \max\limits_{L\in\mathcal{F}}(C_L)\frac{(2C)^p}{r^p}\cdot \max\limits_{L\in\mathcal{F}}\vol(L)\leq \frac{2^pC^{p+2}}{r^p},
\end{equation*}
and
\begin{equation*}
\Cov{Cr}{M}\leq k\cdot \frac{2^pC^{p+2}}{r^p} =: B(r).
\end{equation*}
By Corollary \ref{CorollaryCapEstimation},
\begin{equation*}
\Cap{3Cr}{M'}\geq \frac{1}{(6C)^{n'}\cdot (C)^2} \Cap{\frac{r}{2}}{M'}.
\end{equation*} 
Next, by (\ref{AConstant}),
\begin{equation*}
\Cap{3Cr}{M'}\geq\frac{k}{(6C)^{n'}\cdot C^4} \cdot \frac{2^{p'}}{r^{p'}} =:A(r)
\end{equation*}
It follows that
\begin{equation*}
\frac{A(r)}{B(r)} = \beta\cdot r^{p-p'}.
\end{equation*}
where $\beta$ depends only on $C$, $p'$, and $n'$. Since $p<p'$ then $\lim\limits_{r\to\ 0} \frac{A(r)}{B(r)} = +\infty$. Hence, there exists $\epsilon_0<r$ such that $\Cap{3C\epsilon_0}{M'}\geq A(r)>B(r) \geq \Cov{C\epsilon_0}{M}$. By Lemma \ref{GHlemma-converse}, we obtain $d_{{\rm GH}}(M,M')>C\epsilon_0$.
\end{proof}

\end{document}